\documentclass[12pt]{article}
\usepackage{amsmath,amssymb,amsthm}
\usepackage[margin=1in]{geometry}
\usepackage{hyperref}
\usepackage{tikz}
\definecolor{xcol}{RGB}{214,96,30}
\definecolor{ycol}{RGB}{40,90,160}
\definecolor{nbr}{RGB}{252,226,205}

\newtheorem{thm}{Theorem}[section]
\newtheorem{lem}[thm]{Lemma}
\newtheorem{prop}[thm]{Proposition}

\newtheorem{claim}{Claim}
\theoremstyle{definition}

\title{Graphs with three distinct distance eigenvalues, one of which is zero}
\author{Suil O\thanks{Department of Applied Mathematics and Statistics, The State University of New York, Korea, Incheon, 21985, suil.o@sunykorea.ac.kr. Research supported by the National Research Foundation of Korea (NRF) grant funded by the Korea government(MSIT) No. RS-2025-23523950.}}
\date{}

\begin{document}
\maketitle

\begin{abstract}
Koolen, Hayat and Iqbal asked for the connected graphs whose distance matrix has exactly three distinct eigenvalues, the middle one being $0$. The examples they gave, namely the Johnson graphs, the Hamming graphs and the halved cubes, are distance-regular. In this note, we show that the distance spectrum of such a graph is $\{\rho,0,-c\}$ with integers $c\ge2$ and $c\mid\rho$. For $c=2$ the graphs are exactly the cocktail party graphs, and for $c=3$ they are exactly the Cartesian product $K_3\,\square\,K_3$ of two triangles, the Petersen graph and its complement. We also show that such a graph need not be transmission-regular: Seidel switching the lattice graph $K_6\,\square\,K_6$ with respect to a line and the triangular graph $T(11)$ with respect to a maximal clique of size $10$ gives graphs on $36$ and $55$ vertices with distance spectra $\{54,0^{(26)},(-6)^{(9)}\}$ and $\{81,0^{(45)},(-9)^{(9)}\}$. Among the graphs of diameter $2$ obtained by switching a strongly regular graph with smallest eigenvalue $-2$ with respect to one part of an equitable partition, these two are the only ones that have exactly three distance eigenvalues and are not transmission-regular. We also observe that the transmission-regular graphs with exactly three distinct distance eigenvalues, the smallest one being $-3$, are exactly the complements of the strongly regular graphs with smallest eigenvalue $-2$ other than the cocktail party graphs.


\end{abstract}

\section{Introduction}

All graphs in this note are finite and simple. For a connected graph $G$ with vertices $v_1,\dots,v_n$, the \emph{distance matrix} $D(G)$ is the $n\times n$ matrix whose $(i,j)$ entry is the distance $d(v_i,v_j)$, and its eigenvalues are the \emph{distance eigenvalues} of $G$. The \emph{transmission} of a vertex $u$ is $\sum_{w}d(u,w)$, and $G$ is \emph{transmission-regular} if all vertices have the same transmission, that is, if $\mathbf 1$ is an eigenvector of $D(G)$. For background on distance spectra we refer to~\cite{AH,AABC}, and for standard terminology to~\cite{West,BH,GR}.

The complete graphs are the only connected graphs with exactly two distinct distance eigenvalues (see~\cite{ZL}), so the next case is three. Koolen, Hayat and Iqbal~\cite{KHI} proposed the following problem.

\medskip
\noindent\emph{Problem}~\cite[Problem~6.1]{KHI}. Determine the connected graphs with exactly three distinct distance eigenvalues, the middle one being $0$.
\medskip

They noted that the Johnson graphs, the Hamming graphs and the halved cubes are examples, and that for distance-regular graphs the problem was solved by Aalipour et al.~\cite{AABC}, building on the work of Koolen and Shpectorov~\cite{KS} on distance-regular graphs whose distance matrix has exactly one positive eigenvalue. Zhang and Lin~\cite{ZL} solved a related problem from~\cite{KHI} for graphs without $C_3$ and $C_4$: such a graph has exactly three distinct distance eigenvalues with smallest one $-3$ if and only if it is the Petersen graph, whose distance spectrum is $\{15,0^{(4)},(-3)^{(5)}\}$.

Since $D(G)$ has zero trace and every off-diagonal entry is positive, a graph as in the problem has distance spectrum $\{\rho,0,-c\}$ with $\rho>0>-c$ and $\rho$ simple; we write $a$ and $b$ for the multiplicities of $0$ and $-c$. In this note, we determine these graphs for $c\le 3$. The \emph{lattice graph} $L_2(q)$ is the Cartesian product $K_q\,\square\,K_q$, the graph $T(m)$ is the line graph of $K_m$, called the \emph{triangular graph}, and the \emph{cocktail party graph} $CP(k)$ is $K_{2k}$ minus a perfect matching.

\begin{thm}\label{thm:main}
Let $G$ be a connected graph with distance spectrum $\{\rho,0,-c\}$, where $\rho>0>-c$.
\begin{itemize}
\item[(a)] $c$ is an integer with $c\ge 2$, $b\ge 2$ and $\rho=bc$. If $G$ has diameter at least $3$, then $c\ge4$.
\item[(b)] If $c=2$, then $G$ is the cocktail party graph $CP(k)$ for some $k\ge2$.
\item[(c)] If $c=3$, then $G$ is $L_2(3)$, the Petersen graph or $T(5)$.
\end{itemize}
\end{thm}

Problem~6.2 of~\cite{KHI} asks for the connected graphs with exactly three distinct distance eigenvalues, the smallest one being $-3$. Part (c) answers this question when the middle eigenvalue is $0$. When the graph is transmission-regular, the question is answered by Seidel's classification of the strongly regular graphs with smallest eigenvalue $-2$ (Theorem~\ref{thm:seidel} below).

\begin{prop}\label{prop:tr}
A connected transmission-regular graph $G$ has exactly three distinct distance eigenvalues, the smallest one being $-3$, if and only if $\overline G$ is a strongly regular graph with smallest eigenvalue $-2$ other than a cocktail party graph, that is, $G$ is $\overline{T(m)}$ with $m\ge5$, $\overline{L_2(q)}$ with $q\ge3$, $T(5)$, the Clebsch graph, or the complement of the Shrikhande graph, the Schl\"afli graph or a Chang graph.
\end{prop}

The middle distance eigenvalues of these graphs are $m-5$, $q-3$, $0$, $1$, $1$, $3$ and $3$, respectively, so the middle eigenvalue is $0$ exactly for $T(5)$, the Petersen graph and $L_2(3)$, in accordance with part (c).

\begin{figure}[ht]
\centering
\begin{minipage}[b]{0.47\textwidth}\centering
\begin{tikzpicture}[scale=0.62,
  v/.style={circle,fill=ycol,inner sep=0pt,minimum size=5pt},
  xv/.style={circle,fill=xcol,inner sep=0pt,minimum size=6pt}]
  \fill[nbr,rounded corners=4pt] (0.6,-0.6) rectangle (5.4,-5.4);
  \foreach \i in {0,...,5} {
    \draw[gray!45] (0,-\i) -- (5,-\i);
    \draw[gray!45] (\i,0) -- (\i,-5);
  }
  \draw[xcol!60,line width=2.2pt,rounded corners=3pt,opacity=0.35] (0,0.35) -- (0,-5.35);
  \foreach \i in {0,...,5} {
    \node[xv] at (0,-\i) {};
    \foreach \j in {1,...,5} { \node[v] at (\j,-\i) {}; }
  }
  \draw[xcol,line width=1pt] (0,0) circle (0.27);
  \foreach \j in {1,...,5} { \draw[gray,line width=0.6pt] (\j,0) circle (0.2); }
  \node[left,xcol] at (-0.35,0) {$u$};
  \node[below,xcol] at (0,-5.45) {$X$};
\end{tikzpicture}

\small (a) $G_{36}$
\end{minipage}\hfill
\begin{minipage}[b]{0.47\textwidth}\centering
\begin{tikzpicture}[scale=0.95]
  \def\R{2.2}
  \foreach \a in {1,...,10} { \coordinate (P\a) at ({90-36*(\a-1)}:\R); }
  \foreach \a in {2,...,10} {
    \foreach \b in {\a,...,10} {
      \ifnum\a<\b \draw[ycol!70,line width=0.45pt] (P\a) -- (P\b); \fi
    }
  }
  \foreach \b in {2,...,10} { \draw[gray!60,dashed,line width=0.45pt] (P1) -- (P\b); }
  \foreach \a in {1,...,10} { \node[circle,fill=xcol,inner sep=0pt,minimum size=6pt] at (P\a) {}; }
  \draw[xcol,line width=1pt] (P1) circle (0.13);
  \node[above,xcol] at (P1) {$u$};
\end{tikzpicture}

\small (b) $G_{55}$
\end{minipage}
\caption{(a) $G_{36}$. The $36$ vertices form a $6\times6$ grid, and every row and every column of $L_2(6)$ is a clique (grey lines). The left column is the switching set $X$; it remains a clique. A vertex $u\in X$ is adjacent to the rest of $X$ and to the $25$ vertices in the shaded block, but not to the five circled vertices of its own row. (b) $G_{55}$. Write $T(11)$ on the pairs from $\{p,1,\dots,10\}$. The ten points stand for the clique $X=\{\{p,i\}:1\le i\le10\}$, and the $45$ chords stand for the pairs $\{i,j\}$, which induce $T(10)$. The vertex $\{p,i\}$ is adjacent to the chords avoiding the point $i$. For the circled point $u$ these are the $36$ solid chords; the dashed chords through $u$ are its non-neighbors.}
\label{fig:examples}
\end{figure}
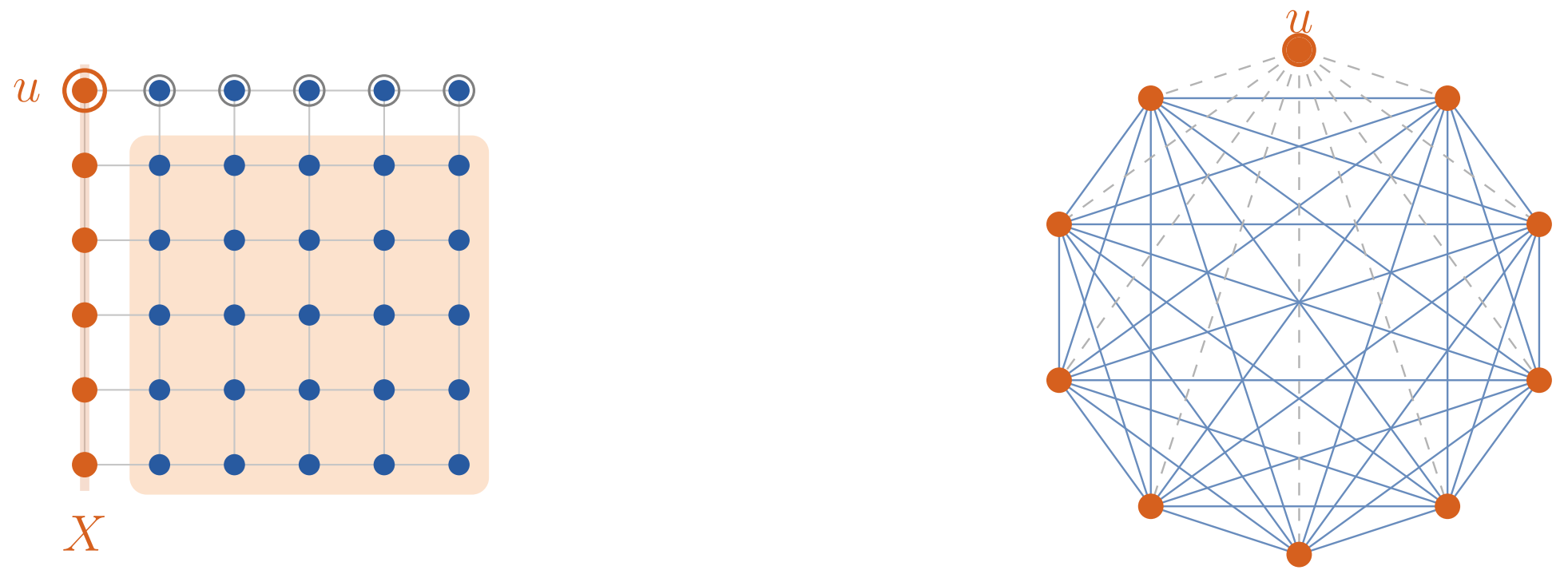

All graphs in (b) and (c) are strongly regular, so they are transmission-regular, as are the examples listed in~\cite{KHI}. A natural guess is that every graph in the problem is transmission-regular. Our second result shows that this is false. For a graph $H$ and $X\subseteq V(H)$, the \emph{Seidel switch} of $H$ with respect to $X$ is obtained by exchanging edges and non-edges between $X$ and $V(H)\setminus X$. Let $G_{36}$ be the Seidel switch of $L_2(6)$ with respect to a line, that is, the six vertices of one column of the $6\times6$ grid, and let $G_{55}$ be the Seidel switch of $T(11)$ with respect to the ten pairs containing a fixed element of $\{1,\dots,11\}$. They are drawn in Figure~\ref{fig:examples}.

\begin{thm}\label{thm:switch}
The graphs $G_{36}$ and $G_{55}$ have diameter $2$ and distance spectra
\[
\operatorname{spec}D(G_{36})=\{54,\ 0^{(26)},\ (-6)^{(9)}\},\qquad
\operatorname{spec}D(G_{55})=\{81,\ 0^{(45)},\ (-9)^{(9)}\},
\]
and neither of them is transmission-regular. Moreover, let $H$ be a connected strongly regular graph with smallest eigenvalue $-2$, let $\emptyset\ne X\subsetneq V(H)$ be such that $\{X,V(H)\setminus X\}$ is an equitable partition of $H$, and let $G$ be the Seidel switch of $H$ with respect to $X$. If $G$ has diameter $2$, has exactly three distinct distance eigenvalues and is not transmission-regular, then $G$ is $G_{36}$ or $G_{55}$.
\end{thm}

The proof of Theorem~\ref{thm:switch} shows a little more: in every other case where such a switch $G$ has exactly three distinct distance eigenvalues, $G$ is transmission-regular. A transmission-regular graph of diameter $2$ is regular, so then $G$ is again a strongly regular graph with smallest eigenvalue $-2$.

The note is organized as follows. Section~2 collects the tools, and Section~3 proves Theorems~\ref{thm:main} and~\ref{thm:switch} and Proposition~\ref{prop:tr}. Section~4 contains concluding remarks.

\section{Definitions and Tools}

We write $J$ for the all-ones matrix, $\mathbf 1$ for the all-ones vector, $\chi_X$ for the characteristic vector of a set $X$, and $e_u$ for the standard basis vector of a vertex $u$. For a symmetric matrix $M$, $\lambda_{\min}(M)$ denotes its smallest eigenvalue. We use Cauchy interlacing in the following form: if $N$ is a principal submatrix of a real symmetric matrix $M$, then $\lambda_{\min}(M)\le\lambda_{\min}(N)$, and $N$ has at most as many positive eigenvalues as $M$.

A partition $\{X,V\setminus X\}$ of the vertex set of a $k$-regular graph $H$ is \emph{equitable} if every vertex of $X$ has the same number $p$ of neighbors in $X$ and every vertex of $V\setminus X$ has the same number $r$ of neighbors in $X$. In this case $\eta=p-r$ is an eigenvalue of $A(H)$ with eigenvector $\chi_X-\frac{|X|}{n}\mathbf 1$, and counting edges between $X$ and $V\setminus X$ gives $r=|X|(k-\eta)/n$~\cite[Section~9.3]{GR}.

The Seidel matrix of $H$ is $S_H=J-I-2A(H)$. If $G$ is the Seidel switch of $H$ with respect to $X$, $\sigma=\mathbf 1-2\chi_X$ and $\Sigma=\operatorname{diag}(\sigma)$, then $S_G=\Sigma S_H\Sigma$. If $G$ has diameter at most $2$, then $D(G)=A(G)+2(J-I-A(G))$, so
\begin{equation}\label{eq:seidel}
2D(G)=S_G-3I+3J.
\end{equation}

We use the following two results.

\begin{thm}[\cite{LHWS}]\label{thm:LHWS}
Let $G$ be a connected graph on $n$ vertices. Then $\lambda_{\min}(D(G))=-2$ with multiplicity $n-k$ if and only if $G$ is a complete $k$-partite graph with $2\le k\le n-1$.
\end{thm}

\begin{thm}[\cite{Seidel}; see also~\cite{BV}]\label{thm:seidel}
A connected strongly regular graph with smallest eigenvalue $-2$ is one of $T(m)$ with $m\ge5$, $L_2(q)$ with $q\ge3$, $CP(k)$ with $k\ge2$, the Petersen graph, the Shrikhande graph, the complement of the Clebsch graph, the Schl\"afli graph, and the three Chang graphs.
\end{thm}

\section{Main Results}

We first prove part (a) of Theorem~\ref{thm:main}.

\begin{lem}\label{lem:basic}
Let $G$ be a connected graph with distance spectrum $\{\rho,0,-c\}$, where $\rho>0>-c$. Then $c$ is an integer, $c\ge2$, $b\ge2$ and $\rho=bc$. If $G$ has diameter at least $3$, then $c\ge 4$.
\end{lem}

\begin{proof}
Since $G$ has three distinct distance eigenvalues, $G$ is not complete and $n\ge3$. Every off-diagonal entry of $D=D(G)$ is positive, so $D^2$ is a positive matrix and $D$ is primitive. By the Perron--Frobenius theorem, $\rho$ is simple and $|\lambda|<\rho$ for every other eigenvalue $\lambda$. As $\operatorname{tr}D=0$, we have $\rho=bc$, and $b=1$ would give the eigenvalue $-\rho$; hence $b\ge2$. The characteristic polynomial of $D$ is $t^a(t-\rho)(t+c)^b$, so $(t-\rho)(t+c)^b$ has integer coefficients, and every automorphism of $\overline{\mathbb Q}$ permutes its roots and preserves their multiplicities. Since $b\ne1$, every such automorphism fixes $-c$, and hence also $\rho=bc$. Hence $\rho$ and $c$ are rational algebraic integers, that is, integers; this also follows from~\cite[Theorem~3.1]{KHI}.

Since $G$ is connected and not complete, it has vertices $u,v,w$ with $u\sim v\sim w$ and $d(u,w)=2$. The principal submatrix of $D$ on $\{u,v,w\}$ has eigenvalue $-2$ with eigenvector $e_u-e_w$, so $c\ge2$ by interlacing. If $G$ has diameter at least $3$, a shortest path on four vertices spans a principal submatrix equal to $D(P_4)$, whose characteristic polynomial is $(t^2-4t-6)(t^2+4t+2)$. Its smallest eigenvalue is $-2-\sqrt2$, so $c\ge2+\sqrt2$, and $c\ge4$ since $c$ is an integer.
\end{proof}

\begin{proof}[Proof of Theorem~\ref{thm:main}(b)]
Let $c=2$, and let $m=b$ be the multiplicity of $-2$, the smallest distance eigenvalue. Since $\rho$ is simple and $0$ is an eigenvalue, $1\le m\le n-2$. By Theorem~\ref{thm:LHWS}, $G=K_{n_1,\dots,n_k}$ with $k=n-m\ge2$. Let $B$ be the block-diagonal matrix with all-ones blocks on the parts. Then $D+2I=J+B$, which vanishes on every vector summing to $0$ on each part. On the span of the normalized characteristic vectors of the parts, $J+B$ acts as the $k\times k$ matrix $M=\operatorname{diag}(n_1,\dots,n_k)+zz^{\top}$, where $z_i=\sqrt{n_i}$. The matrix $M$ is positive definite, and its eigenvalues are distance eigenvalues of $G$ shifted by $2$. Hence the eigenvalues of $M$ lie in $\{\rho+2,2\}$ and $\rho+2$ is simple, so $M-2I$ has rank $1$. Therefore every $2\times2$ principal minor of $M-2I$ vanishes:
\[
(2n_i-2)(2n_j-2)=n_in_j,\qquad\text{that is,}\qquad h(n_i)h(n_j)=1\ \text{ with } h(t)=2-\tfrac2t .
\]
Since $h(1)=0$, $h(2)=1$ and $1<h(t)<2$ for $t\ge3$, all $n_i$ equal $2$, and $G=CP(k)$. Conversely, $CP(k)$ has distance spectrum $\{2k,0^{(k-1)},(-2)^{(k)}\}$.
\end{proof}

For the case $c=3$, let $M(F)=2(J-I)-A(F)$ for a graph $F$, and let $\mathcal H$ be the class of graphs $F$ such that $\lambda_{\min}(M(F))\ge-3$ and $M(F)$ has at most one positive eigenvalue. By interlacing, $\mathcal H$ is closed under taking induced subgraphs. Write $\mathcal H(n)$ for the set of graphs in $\mathcal H$ on $n$ vertices. For $t,r\ge0$ let $M_{t,r}=tK_2\cup rK_1$, and let $\mathcal F(n)$ be the set of graphs on $n$ vertices of the following four types:
\begin{itemize}
\item[(i)] the complement of $M_{t,n-2t}$;
\item[(ii)] $K_s\vee M_{t,n-s-2t}$ with $s\in\{1,2\}$;
\item[(iii)] $M_{t,n-2t}$;
\item[(iv)] $(K_1\vee M_{t,n-2-2t})\cup K_1$.
\end{itemize}
We call types (ii), (iii) and (iv) \emph{sparse}. The \emph{apexes} of a sparse graph are the vertices of $K_s$ in type (ii) and the vertex of $K_1$ joined to $M_{t,n-2-2t}$ in type (iv); type (iii) has no apex. In type (i) every vertex has degree at least $n-2$. In a sparse graph, every apex has degree at least $n-2$, every other vertex has degree at most $3$, and deleting the apexes leaves a graph of maximum degree at most $1$.

\begin{proof}[Proof of Theorem~\ref{thm:main}(c)]
Let $c=3$. By Lemma~\ref{lem:basic}, $G$ has diameter $2$, so $D(G)[Y]=M(G[Y])$ for every $Y\subseteq V(G)$. Since $D(G)$ has smallest eigenvalue $-3$ and exactly one positive eigenvalue, interlacing gives $G[Y]\in\mathcal H$ for every $Y$; in particular $G\in\mathcal H(n)$.

\begin{claim}\label{cl:comp}
$\mathcal H(19)=\mathcal F(19)$, and the connected graphs on at most $19$ vertices with distance spectrum $\{\rho,0,-3\}$ are exactly $L_2(3)$, the Petersen graph and $T(5)$.
\end{claim}

\begin{proof}
This is a computer search. Since $\mathcal H$ is closed under taking induced subgraphs, every graph in $\mathcal H(n+1)$ arises from a graph in $\mathcal H(n)$ by adding one vertex. We generated $\mathcal H(n)$ for $n\le19$ in this way, removing isomorphic copies with nauty~\cite{MP}. Membership was decided exactly from the integer characteristic polynomials: since $M(F)$ is symmetric, its characteristic polynomial has only real roots, so Descartes' rule of signs counts its positive eigenvalues exactly, and $M(F)+3I$ is positive semidefinite if and only if the coefficients of its characteristic polynomial alternate in sign. Before the exact test, a candidate was discarded if a floating-point computation placed its smallest eigenvalue below $-3-10^{-6}$ or its second largest eigenvalue above $10^{-6}$; for integer matrices of order at most $19$ the floating-point error is far below $10^{-6}$, so no graph of $\mathcal H$ is lost in this step. No candidate passing this filter failed the exact test. The numbers of graphs are
\[
|\mathcal H(n)|=2,4,11,28,78,191,391,639,770,674,465,276,135,76,53,43,47,48
\]
for $n=2,\dots,19$. The $48$ graphs in $\mathcal H(19)$ are exactly the $48$ graphs in $\mathcal F(19)$.\footnote{The code and the list of the $3931$ graphs in $\bigcup_{n\le19}\mathcal H(n)$ are available from the author upon request.} Among all graphs found, the connected ones with distance spectrum $\{\rho,0,-3\}$ are $L_2(3)$, the Petersen graph and $T(5)$, with spectra $\{12,0^{(4)},(-3)^{(4)}\}$, $\{15,0^{(4)},(-3)^{(5)}\}$ and $\{12,0^{(5)},(-3)^{(4)}\}$.
\end{proof}

\begin{claim}\label{cl:induction}
$\mathcal H(n)\subseteq\mathcal F(n)$ for every $n\ge19$.
\end{claim}

\begin{proof}
We use induction on $n$; the case $n=19$ is Claim~\ref{cl:comp}. Let $n=N+1\ge20$ and $F\in\mathcal H(n)$. By the induction hypothesis, $F-v\in\mathcal F(N)$ for every vertex $v$.

Suppose first that $F-v$ is of type (i) for some $v$. For every $w\ne v$, each of the $N-1\ge18$ vertices of $V\setminus\{v,w\}$ has degree at least $N-2$ in $F-v$, hence degree at least $N-3\ge16$ in $F-w$. A sparse graph has at most two vertices of degree at least $4$, so every $F-w$ is of type (i), that is, every $\overline F-w$ is a matching. If a vertex $u$ had two neighbors $u_1,u_2$ in $\overline F$, then for any $w\notin\{u,u_1,u_2\}$ the graph $\overline F-w$ would not be a matching. Hence $\overline F$ is a matching, and $F$ is of type (i).

Now suppose that every $F-v$ is sparse. Let $B=\{u:\deg_F(u)\ge N-2\}$. Since $\deg_{F-v}(u)\in\{\deg_F(u),\deg_F(u)-1\}$, every apex of $F-v$ lies in $B$, and every non-apex $u$ of $F-v$ has $\deg_F(u)\le4<N-2$. So the set of apexes of $F-v$ is $B\setminus\{v\}$. If $B=V$, then $F-v$ would have $N\ge 3$ apexes, which is impossible. Choosing $v\notin B$, we get $|B|\le2$. For $v\notin B$, the graph $(F-B)-v$ is $F-v$ with its apexes deleted, so it has maximum degree at most $1$; since $|V\setminus B|\ge18$, the graph $F-B$ has maximum degree at most $1$. If $B=\emptyset$, then $F$ is of type (iii). If $B=\{b_1,b_2\}$, then every $F-v$ with $v\notin B$ is of type (ii) with $s=2$; using two such vertices $v$, we see that $b_1b_2$ is an edge and both $b_1,b_2$ are adjacent to all of $V\setminus B$, so $F=K_2\vee(F-B)$ is of type (ii). Finally, let $B=\{u_0\}$, and let $Z$ be the set of non-neighbors of $u_0$ in $V\setminus B$. Every $F-v$ with $v\notin B$ is of type (ii) with $s=1$ or of type (iv), so $|Z\setminus\{v\}|\le1$; taking $v\notin Z\cup B$ gives $|Z|\le1$. If $Z=\emptyset$, then $F=K_1\vee(F-u_0)$ is of type (ii). If $Z=\{w\}$, then every $F-v$ with $v\notin\{u_0,w\}$ is of type (iv) with $w$ as its isolated vertex, so $w$ is isolated in $F$ and $F=(K_1\vee(F-u_0-w))\cup K_1$ is of type (iv).
\end{proof}

\begin{claim}\label{cl:fam}
No connected graph in $\mathcal F(n)$ has distance spectrum $\{\rho,0,-3\}$.
\end{claim}

\begin{proof}
Graphs of types (iii) and (iv) are disconnected. A graph of type (i) is complete if $t=0$, and then it has only two distance eigenvalues; otherwise it is complete multipartite but not complete, so its smallest distance eigenvalue is $-2$ by Theorem~\ref{thm:LHWS}. Let $F$ be of type (ii). If $aa'$ is an edge of the matching $tK_2$, then every other vertex is at the same distance from $a$ and $a'$, so $D(F)(e_a-e_{a'})=-(e_a-e_{a'})$. The same holds for the two apexes if $s=2$. Since $-1$ is not in the spectrum, $t=0$ and $s=1$, so $F=K_{1,n-1}$, which is complete bipartite and has smallest distance eigenvalue $-2$.
\end{proof}

If $n\le19$, Claim~\ref{cl:comp} gives the result. If $n\ge20$, then $G\in\mathcal F(n)$ by Claim~\ref{cl:induction}, which contradicts Claim~\ref{cl:fam}.
\end{proof}

We next prove Proposition~\ref{prop:tr}.

\begin{proof}[Proof of Proposition~\ref{prop:tr}]
Let $G$ be a connected transmission-regular graph with exactly three distinct distance eigenvalues, the smallest one being $-3$. If $G$ had diameter at least $3$, then, as in the proof of Lemma~\ref{lem:basic}, interlacing with $D(P_4)$ would give a distance eigenvalue at most $-2-\sqrt2<-3$. So $G$ has diameter $2$, and the transmission of a vertex is $2(n-1)$ minus its degree; hence $G$ is $k$-regular for some $k$. Then $D(G)=2J-2I-A(G)$ has eigenvalue $2n-2-k$ on $\mathbf 1$ and eigenvalue $-2-\lambda$ on every eigenvector of $A(G)$ orthogonal to $\mathbf 1$ with eigenvalue $\lambda$. Since $-2-\lambda=2n-2-k$ would give $\lambda=k-2n<-k$, the adjacency matrix $A(G)$ has exactly three distinct eigenvalues $k>\theta>\tau$, so $G$ is strongly regular, and the smallest distance eigenvalue is $-2-\theta$. Hence $\theta=1$. The complement $\overline G$ is regular with eigenvalues $n-1-k$, $-1-\tau$ and $-1-\theta=-2$, and it is connected, since otherwise $G$ would be complete multipartite and would have $\theta=0$. So $\overline G$ is one of the graphs in Theorem~\ref{thm:seidel}, and it is not a cocktail party graph, because $\overline{CP(k)}=kK_2$ is disconnected.

Conversely, let $\overline G$ be a graph in Theorem~\ref{thm:seidel} other than a cocktail party graph. The cocktail party graphs are the only complete multipartite graphs in the list, so $G$ is connected, and $G$ is strongly regular with eigenvalues $n-1-k$, $\theta=1$ and $\tau<1$. Thus $G$ has diameter $2$ and distance eigenvalues $2n-2-k$, $-3$ and $-2-\tau>-3$, which are three distinct values. The middle eigenvalues listed after Proposition~\ref{prop:tr} are the values $-2-\tau=-1+\theta(\overline G)$, read off from $\theta(\overline G)=m-4$, $q-2$, $1$, $2$, $2$, $4$ and $4$.
\end{proof}

We now turn to Seidel switching. The next lemma computes the distance spectrum of the switch in terms of a $2\times2$ matrix.

\begin{lem}\label{lem:switch}
Let $H$ be a connected strongly regular graph on $n$ vertices with valency $k$ and eigenvalues $k>\theta>-2$, of multiplicities $1$, $m_\theta$ and $m_{-2}$. Put $c=\theta+2$ and $\lambda_1=n-4-2k$. Let $\emptyset\ne X\subsetneq V(H)$ be such that $\{X,V(H)\setminus X\}$ is equitable, let $\eta\in\{\theta,-2\}$ be the eigenvalue of $A(H)$ with eigenvector $y=\chi_X-\frac{|X|}{n}\mathbf 1$, and put $\varepsilon=-4-2\eta$, $s=|X|$ and
\[
Y_U=\begin{pmatrix}\lambda_1&0\\0&\varepsilon\end{pmatrix}+3\,zz^{\top},\qquad z=\Bigl(\frac{n-2s}{\sqrt n},\ -2\sqrt{\frac{s(n-s)}{n}}\Bigr)^{\!\top}.
\]
Let $G$ be the Seidel switch of $H$ with respect to $X$, and suppose that $G$ has diameter at most $2$. Then the spectrum of $2D(G)$ consists of the two eigenvalues of $Y_U$, the eigenvalue $-2c$ with multiplicity $m_\theta-1$ if $\eta=\theta$ and $m_\theta$ otherwise, and the eigenvalue $0$ with multiplicity $m_{-2}-1$ if $\eta=-2$ and $m_{-2}$ otherwise. Moreover, $G$ is transmission-regular if and only if $s=n/2$ or $\lambda_1=\varepsilon$.
\end{lem}

\begin{proof}
Let $\sigma=\mathbf 1-2\chi_X$ and $\Sigma=\operatorname{diag}(\sigma)$. Since $\Sigma J\Sigma=\sigma\sigma^{\top}$ and $S_G=\Sigma S_H\Sigma$, equation~\eqref{eq:seidel} gives $2D(G)=\Sigma Y\Sigma$ with $Y=S_H-3I+3\sigma\sigma^{\top}$. As $\Sigma^2=I$, the matrices $2D(G)$ and $Y$ have the same spectrum. The matrix $S_H-3I=J-4I-2A(H)$ has eigenvalue $n-4-2k=\lambda_1$ on $\mathbf 1$, eigenvalue $-4-2\theta=-2c$ on the $\theta$-eigenspace $V_\theta$, and eigenvalue $0$ on the $(-2)$-eigenspace $V_{-2}$; in particular it has eigenvalue $\varepsilon$ on $y$. Let $U=\operatorname{span}\{\mathbf 1,y\}$. Since $\sigma=(1-\frac{2s}{n})\mathbf 1-2y\in U$, the matrix $Y$ agrees with $S_H-3I$ on $U^{\perp}$, which gives the eigenvalues $-2c$ and $0$ with the stated multiplicities. We have $\|y\|^2=s(n-s)/n$, so in the orthonormal basis $\mathbf 1/\sqrt n$, $y/\|y\|$ of $U$ the vector $\sigma$ has coordinates $z$, and $Y$ acts on $U$ as $Y_U$.

Since $\Sigma\mathbf 1=\sigma$ and $\Sigma\sigma=\mathbf 1$, the vector $\mathbf 1$ is an eigenvector of $D(G)$ if and only if $\sigma$ is an eigenvector of $Y$, that is, if and only if $z$ is an eigenvector of $Y_U$. As $\|z\|^2=n$, we have $Y_Uz=\operatorname{diag}(\lambda_1,\varepsilon)z+3nz$, so this happens if and only if $z$ is an eigenvector of $\operatorname{diag}(\lambda_1,\varepsilon)$. Since $z_2\ne0$, this means $z_1=0$ or $\lambda_1=\varepsilon$, that is, $s=n/2$ or $\lambda_1=\varepsilon$.
\end{proof}

\begin{proof}[Proof of Theorem~\ref{thm:switch}]
We first treat $G_{36}$ and $G_{55}$. In $L_2(q)$, the parameters are $n=q^2$, $k=2q-2$ and $\theta=q-2$. A vertex outside a line $X$ has exactly one neighbor in $X$, and a vertex of $X$ has $q-1$ neighbors in $X$, so $\{X,V\setminus X\}$ is equitable with $\eta=q-2=\theta$. In $T(m)$, the parameters are $n=\binom m2$, $k=2m-4$ and $\theta=m-4$. Let $X$ be the set of the $m-1$ pairs containing a fixed element $p$. A pair not containing $p$ has exactly two neighbors in $X$, and a pair in $X$ has $m-2$ neighbors in $X$, so $\{X,V\setminus X\}$ is equitable with $\eta=m-4=\theta$.

For $H=L_2(6)$ we have $n=36$, $k=10$, $c=6$, $\lambda_1=12$, $s=6$, $m_\theta=10$ and $m_{-2}=25$. For $H=T(11)$ we have $n=55$, $k=18$, $c=9$, $\lambda_1=15$, $s=10$, $m_\theta=10$ and $m_{-2}=44$. In both cases $\eta=\theta$, so $\varepsilon=-2c$, and
\begin{equation}\label{eq:Delta}
\det Y_U=-2c\lambda_1+\frac3n\bigl(4\lambda_1s(n-s)-2c(n-2s)^2\bigr)
\end{equation}
equals $-144+\frac{3}{36}(4\cdot12\cdot180-12\cdot576)=0$ and $-270+\frac{3}{55}(4\cdot15\cdot450-18\cdot1225)=0$, respectively. Hence the eigenvalues of $Y_U$ are $0$ and $\operatorname{tr}Y_U=\lambda_1-2c+3n$, which equals $108$ and $162$. Before applying Lemma~\ref{lem:switch} we check that the switches have diameter $2$.

In $G_{36}$ the column $X$ is still a clique, and a vertex $u\in X$ is adjacent exactly to the other vertices of $X$ and to the $25$ vertices outside $X$ and outside the row of $u$. The edges outside $X$ are those of $L_2(6)$, so two non-adjacent vertices outside $X$ lie in different rows, and the vertex of $X$ in a third row is adjacent to both. A vertex $v\notin X$ and the vertex of $X$ in its row have any other vertex of $X$ as a common neighbor, and all other pairs of vertices are adjacent. So $G_{36}$ has diameter $2$, and its degrees are $14$ outside $X$ and $30$ on $X$. In $G_{55}$ the set $X=\{\{p,i\}:i\ne p\}$ is still a clique, and $\{p,i\}$ is adjacent exactly to the other vertices of $X$ and to the pairs not containing $p$ or $i$. Two disjoint pairs $\{i,j\},\{k,l\}$ have the common neighbor $\{p,r\}$ for any $r\notin\{p,i,j,k,l\}$, and $\{i,j\}$ and $\{p,i\}$ have the common neighbor $\{p,r\}$ for any $r\notin\{p,i,j\}$. So $G_{55}$ has diameter $2$, and its degrees are $24$ outside $X$ and $45$ on $X$.

By Lemma~\ref{lem:switch}, the spectrum of $2D(G_{36})$ is $\{108,0^{(26)},(-12)^{(9)}\}$ and the spectrum of $2D(G_{55})$ is $\{162,0^{(45)},(-18)^{(9)}\}$. In diameter $2$ the transmission of a vertex is $2(n-1)$ minus its degree, so neither graph is transmission-regular.

Now let $H$, $X$ and $G$ be as in the second part of the theorem, with the notation of Lemma~\ref{lem:switch}, and suppose that $D(G)$ has exactly three distinct eigenvalues and $G$ is not transmission-regular. By Lemma~\ref{lem:switch}, $s\ne n/2$ and $\lambda_1\ne\varepsilon$. If $H=C_4=CP(2)$, then every equitable partition with $\eta=-2$ consists of two pairs of antipodal vertices, and the switch is $4K_1$, which is disconnected; every equitable partition with $\eta=\theta$ has $s=2=n/2$. So $H\ne C_4$, and $m_\theta,m_{-2}\ge2$ by Theorem~\ref{thm:seidel}. Hence $-2c$ and $0$ are eigenvalues of $2D(G)$. The two eigenvalues $\mu_1,\mu_2$ of $Y_U$ are distinct, because $z_1z_2\ne0$, so $Y_U$ is not a scalar matrix. They are not $0$ and $-2c$, since $\operatorname{tr}D(G)=0$ forces a positive eigenvalue. So exactly one of $\mu_1,\mu_2$ lies in $\{0,-2c\}$, that is, $\det Y_U=0$ or $\det(Y_U+2cI)=0$.

The following table lists $n$, $k$, $\theta$, $\lambda_1$ and $\lambda_1+2c$ for the graphs in Theorem~\ref{thm:seidel}.
\[
\begin{array}{lccccc}
H & n & k & \theta & \lambda_1 & \lambda_1+2c\\ \hline
T(m) & \binom m2 & 2m-4 & m-4 & \frac12(m-1)(m-8) & \frac12m(m-5)\\
L_2(q) & q^2 & 2q-2 & q-2 & q(q-4) & q(q-2)\\
CP(k) & 2k & 2k-2 & 0 & -2k & 4-2k\\
\text{Petersen} & 10 & 3 & 1 & 0 & 6\\
\text{Shrikhande} & 16 & 6 & 2 & 0 & 8\\
\text{Clebsch complement} & 16 & 10 & 2 & -8 & 0\\
\text{Schl\"afli} & 27 & 16 & 4 & -9 & 3\\
\text{Chang} & 28 & 12 & 4 & 0 & 12
\end{array}
\]

Suppose first that $\eta=-2$, so $\varepsilon=0$. Then $\det Y_U=3\lambda_1z_2^2\ne0$, because $\lambda_1\ne\varepsilon=0$. So $\det(Y_U+2cI)=0$, where
\[
\det(Y_U+2cI)=2c(\lambda_1+2c)+3(\lambda_1+2c)z_2^2+6cz_1^2 .
\]
Since $z_1\ne0$, this is positive whenever $\lambda_1+2c\ge0$, so by the table $H=CP(k)$ with $k\ge3$. The $(-2)$-eigenspace of $CP(k)$ consists of the vectors that are constant on each antipodal pair and sum to $0$. Since $\chi_X\in\operatorname{span}\{\mathbf 1,y\}$, the set $X$ is a union of antipodal pairs, and then $G$ is the disjoint union of $CP(s/2)$ and $CP(k-s/2)$, which is disconnected. This is a contradiction.

Hence $\eta=\theta$ and $\varepsilon=-2c$. Then $\det(Y_U+2cI)=3(\lambda_1+2c)z_2^2\ne0$, because $\lambda_1\ne\varepsilon$. So $\det Y_U=0$, where $\det Y_U$ is given by~\eqref{eq:Delta}. If $\lambda_1=0$, then $\det Y_U=-\frac{6c}{n}(n-2s)^2\ne0$, so by the table $H$ is not the Petersen graph, the Shrikhande graph, a Chang graph, $L_2(4)$ or $T(8)$. Also $H$ is neither the complement of the Clebsch graph nor $T(5)$, because these have $\lambda_1=-2c=\varepsilon$. We treat the remaining graphs. Note that $\chi_X=\frac sn\mathbf 1+y$ lies in $\operatorname{span}\{\mathbf 1\}\oplus V_\theta$.

\emph{The Schl\"afli graph.} Here $r=s(k-\theta)/n=4s/9$ is an integer, so $s\in\{9,18\}$. By~\eqref{eq:Delta} with $n=27$, $c=6$ and $\lambda_1=-9$, both values give $\det Y_U=108-756=-648\ne0$.

\emph{$CP(k)$ with $k\ge3$.} Here $V_\theta=V_0$ consists of the vectors $f$ with $f(u')=-f(u)$ for antipodal vertices $u,u'$. So $\chi_X=t\mathbf 1+f$ with $f\in V_0$. If $t+f(u)=t-f(u)$ for some $u$, then $f(u)=0$ and $t\in\{0,1\}$, and then $t\pm f(w)\in\{0,1\}$ forces $f=0$, so $X$ is empty or everything. Hence $t+f(u)$ and $t-f(u)$ are $0$ and $1$ in some order for every $u$, so $X$ contains exactly one vertex of each antipodal pair, and $s=n/2$, which is excluded.

\emph{$L_2(q)$ with $q\ge3$ and $q\ne4$.} Here $\operatorname{span}\{\mathbf 1\}\oplus V_\theta$ consists of the vectors $(i,j)\mapsto a_i+b_j$. If $a_i\ne a_{i'}$ for some $i,i'$, then for every $j$ the numbers $a_i+b_j$ and $a_{i'}+b_j$ are $0$ and $1$ in some order, so $b_j=-\min\{a_i,a_{i'}\}$ does not depend on $j$. Hence $X$ is a union of $\ell$ rows or a union of $\ell$ columns, with $1\le \ell\le q-1$, and $s=\ell q$. Substituting into~\eqref{eq:Delta} gives
\[
\det Y_U=4q\bigl(3(q-2)\,\ell(q-\ell)-2q(q-1)\bigr).
\]
As $\ell(q-\ell)\ge q-1$, a zero requires $3(q-2)\le 2q$, that is, $q\le6$. For $q=3,5,6$ the equation $3(q-2)\ell(q-\ell)=2q(q-1)$ becomes $\ell(3-\ell)=4$, $9\ell(5-\ell)=40$ and $\ell(6-\ell)=5$, and only the last has solutions, namely $\ell\in\{1,5\}$. Replacing $X$ by $V\setminus X$ does not change $G$, so $H=L_2(6)$, $X$ is a line, and $G=G_{36}$.

\emph{$T(m)$ with $m\ge6$ and $m\ne8$.} Here $\operatorname{span}\{\mathbf 1\}\oplus V_\theta$ consists of the vectors $\{i,j\}\mapsto g(i)+g(j)$, so all sums $g(i)+g(j)$ with $i\ne j$ lie in $\{0,1\}$. If $g$ took three distinct values $g_1<g_2<g_3$, then $g_1+g_2<g_1+g_3<g_2+g_3$ would be three distinct elements of $\{0,1\}$. If $g$ took two distinct values $g_1<g_2$, each at least twice, then $2g_1<g_1+g_2<2g_2$ would be three distinct elements of $\{0,1\}$. If $g$ is constant, then $X$ is empty or everything. Hence $g$ takes a value $g_0$ on all elements except one element $p$. Then $2g_0\in\{0,1\}$ and $g(p)+g_0$ is the other element of $\{0,1\}$. If $g_0=0$, then $g(p)=1$ and $X$ is the set of the $m-1$ pairs containing $p$; if $g_0=\frac12$, then $g(p)=-\frac12$ and $X$ is its complement. Replacing $X$ by $V\setminus X$ does not change $G$, so we may take $s=m-1$, and~\eqref{eq:Delta} gives
\[
\det Y_U=2(m-1)(m-2)(m-11).
\]
This vanishes only for $m=11$, so $H=T(11)$ and $G=G_{55}$.
\end{proof}

\section{Concluding Remarks}

Theorem~\ref{thm:main} settles the problem of Koolen, Hayat and Iqbal for $c\le3$, and in this range every example is strongly regular. For $c\ge4$ the picture changes in two ways. First, by Theorem~\ref{thm:switch} the graphs need not be transmission-regular. Second, graphs of larger diameter appear: the hypercubes, the Johnson graphs, the Hamming graphs, the halved cubes and the Doob graphs are distance-regular examples of diameter at least $3$~\cite{AABC,KHI}, and by Lemma~\ref{lem:basic} they all have $c\ge4$. Since a geodesic of length $d$ spans a principal submatrix $D(P_{d+1})$, the diameter of $G$ is bounded in terms of $c$, but we do not know how to control the graphs of diameter at least $3$ that are not distance-regular.

In diameter $2$, the relation $2D(G)=S_G-3I+3J$ and the fact that $D(G)$ has only one positive eigenvalue imply that every graph in the switching class of $G$ has at most one adjacency eigenvalue below $-2$. This is how $G_{36}$ and $G_{55}$ were found. The situation resembles the work of Cheng, Greaves and Koolen~\cite{CGK} on graphs with three adjacency eigenvalues and second largest eigenvalue at most $1$, where the complement has smallest eigenvalue at least $-2$. The relation above brings the theory of graphs with smallest eigenvalue $-2$ to bear on the problem. The difficulty is that the switching class need not contain a strongly regular graph, and that even when it does, the switching set need not be one part of an equitable partition. We ask whether $G_{36}$ and $G_{55}$ are the only connected graphs with exactly three distinct distance eigenvalues, the middle one being $0$, that are not transmission-regular.

\section*{Declaration of generative AI use}
During the preparation of this work the author used Claude (Anthropic) for language refinement, technical editing, computational verification of examples, and symbolic computations. The author reviewed and edited the content and takes full responsibility for it.


\begin{thebibliography}{99}

\bibitem{AABC} G.~Aalipour, A.~Abiad, Z.~Berikkyzy, J.~Cummings, J.~De~Silva, W.~Gao, K.~Heysse, L.~Hogben, F.H.J.~Kenter, J.C.-H.~Lin, M.~Tait, On the distance spectra of graphs, Linear Algebra Appl. 497 (2016) 66--87.

\bibitem{AH} M.~Aouchiche, P.~Hansen, Distance spectra of graphs: a survey, Linear Algebra Appl. 458 (2014) 301--386.

\bibitem{BH} A.E.~Brouwer, W.H.~Haemers, Spectra of Graphs, Springer, New York, 2012.

\bibitem{BV} A.E.~Brouwer, H.~Van~Maldeghem, Strongly Regular Graphs, Cambridge University Press, Cambridge, 2022.

\bibitem{CGK} X.-M.~Cheng, G.R.W.~Greaves, J.H.~Koolen, Graphs with three eigenvalues and second largest eigenvalue at most $1$, J. Combin. Theory Ser. B 129 (2018) 55--78.

\bibitem{GR} C.~Godsil, G.~Royle, Algebraic Graph Theory, Springer, New York, 2001.

\bibitem{KHI} J.H.~Koolen, S.~Hayat, Q.~Iqbal, Hypercubes are determined by their distance spectra, Linear Algebra Appl. 505 (2016) 97--108.

\bibitem{KS} J.H.~Koolen, S.V.~Shpectorov, Distance-regular graphs the distance matrix of which has only one positive eigenvalue, European J. Combin. 15 (1994) 269--275.

\bibitem{LHWS} H.~Lin, Y.~Hong, J.~Wang, J.~Shu, On the distance spectrum of graphs, Linear Algebra Appl. 439 (2013) 1662--1669.

\bibitem{MP} B.D.~McKay, A.~Piperno, Practical graph isomorphism, II, J. Symbolic Comput. 60 (2014) 94--112.

\bibitem{Seidel} J.J.~Seidel, Strongly regular graphs with $(-1,1,0)$ adjacency matrix having eigenvalue $3$, Linear Algebra Appl. 1 (1968) 281--298.

\bibitem{West} D.B.~West, Introduction to Graph Theory, 2nd ed., Prentice Hall, Upper Saddle River, 2001.

\bibitem{ZL} Y.~Zhang, H.~Lin, Graphs with three distinct distance eigenvalues, Appl. Math. Comput. 445 (2023) 127848.

\end{thebibliography}
\end{document}